\documentclass [12pt]{article}
\usepackage{graphicx,amssymb,amsfonts,latexsym,amsmath,amsthm,times}
\usepackage{epsfig}
\usepackage{color}
\usepackage[english]{babel}
\usepackage[a4paper]{geometry}
\def\lanbox{\hbox{$\, \vrule height 0.25cm width 0.25cm depth 0.01cm \,$}}

\numberwithin{equation}{section}

\begin{document}

\vspace*{1.4cm}

\normalsize \centerline{\Large \bf EXISTENCE OF SOLUTIONS
FOR A CLASS OF}

\medskip

\centerline{\Large\bf INTEGRO-DIFFERENTIAL EQUATIONS WITH THE}

\medskip

\centerline{\Large\bf LOGARITHMIC LAPLACIAN AND TRANSPORT}

\vspace*{1cm}

\centerline{\bf Yuming Chen$^{1}$, Vitali Vougalter$^{2 \ *}$}

\vspace*{0.5cm}

\centerline{$^1$ Department of Mathematics, Wilfrid Laurier University}

\centerline{Waterloo, Ontario, N2L 3C5, Canada}

\centerline{ e-mail: ychen@wlu.ca}

\medskip

\centerline{$^{2 \ *}$ Department of Mathematics, University
of Toronto}

\centerline{Toronto, Ontario, M5S 2E4, Canada}

\centerline{ e-mail: vitali@math.toronto.edu}

\medskip

%*******************************************************************
%ABSTRACT
%*******************************************************************

\vspace*{0.25cm}

\noindent {\bf Abstract:}
In this paper, we consider an integro-differential equation in $L^2(\mathbb{R})$, which  involves the logarithmic Laplacian in the presence of a drift term. The linear operator associated with the problem has the Fredholm property. By using  a fixed point technique, we establish the existence of solutions.

\vspace*{0.25cm}



\vspace*{0.5cm}

\bigskip

\bigskip

%%%%%%%%%%%%%%%%%%%%%%%%%%%%%%%%%%%%%%%%%%%%%%%%%%

\setcounter{section}{1}

\centerline{\bf 1. Introduction}

\medskip

\noindent
The present work deals with the solvability of the following
integro-differential equation,
\begin{equation}
\label{h}
\Big[-\frac{1}{2}\hbox{ln}\Big(-\frac{d^{2}}{dx^{2}}\Big)\Big]u+
b\frac{du}{dx}+au+f(x)+
\int_{-\infty}^{\infty}K(x-y)g(u(y))dy=0, 
\end{equation}
where $x\in \mathbb{R}$ and $a, b\in\mathbb{R}, \ b\neq 0$ are constants. Note that the existence of solutions of  a problem similar to~\eqref{h}  with
the fractional Laplacian in the diffusion term has been discussed by Vougalter and Volpert~\cite{VV21}.

Biological processes  
with the nonlocal consumption of resources and  intra-specific competition can be described by nonlocal reaction-diffusion equations (see, for example, \cite{ABVV10,BNPR09,GVA06,VV130} and
references therein).
In such time-dependent problems the space variable $x$ corresponds to the cell
genotype and $u(x,t)$ denotes the cell density as a function of their genotype
and time.
The evolution of the cell density occurs via the
cell proliferation, mutations, transport, and cell influx/efflux.
The diffusion terms correspond to the change of genotype
due to small random mutations  and the integral term describes large mutations.
The function $g(u)$ denotes the rate of cell birth depending on $u$
(density-dependent proliferation). The kernel $K(x-y)$ gives
the proportion of newly born cells changing their genotype from $y$ to $x$.
It is assumed that it depends on the distance between the genotypes.
The source term $f(x)$ designates
the influx/efflux of cells for different genotypes.
Because the models of this kind describe the distribution of a population
density with respect to the genotype, the existence of stationary solutions
of the nonlocal reaction-diffusion equations
corresponds to the existence of biological species.

Let us introduce the concept of   non-Fredholm operators. Consider the
equation
\begin{equation}
\label{eq1}
 -\Delta u + V(x) u - a u=f,
\end{equation}
where $u \in E= H^{2}({\mathbb R}^{d})$ and  $f \in F=
L^{2}({\mathbb R}^{d})$, $d\in {\mathbb N}$, $a$ is a constant, and
the scalar potential function $V(x)$ either vanishes identically
or tends to $0$ at infinity. If $a \geq 0$, the origin is contained in the
essential spectrum of the operator $A : E \to F$ corresponding to the left
side of problem~\eqref{eq1}. As a consequence, this operator does not
satisfy the Fredholm property. Its image is not closed, and, for $d>1$,
the dimension of its kernel and the codimension of its image are
not finite.
Note that for the Fredholm property we consider bounded but not
$H^{2}({\mathbb R}^{d})$ solutions of the corresponding homogeneous adjoint
equation.
If $V(x)\equiv 0$, problem~\eqref{eq1} has
constant coefficients and we can use the Fourier transform to solve it
explicitly. If $f\in L^{2}({\mathbb R}^{d})$ and $xf\in L^{1}({\mathbb R}^{d})$,
then it admits a unique solution in $H^{2}({\mathbb R}^{d})$ if and only if
\begin{equation}
\label{ocb}
\Bigg(f(x), \frac{e^{ipx}}{(2\pi)^{\frac{d}{2}}}\Bigg)_{L^{2}({\mathbb R}^{d})}=0, \quad
p\in S_{\sqrt{a}}^{d} \quad a.e.
\end{equation}
(see Lemmas 5 and 6 of ~\cite{VV14}). Here $S_{\sqrt{a}}^{d}$ denotes the
sphere in ${\mathbb R}^{d}$ of radius $\sqrt{a}$ centered at the origin.
Thus  though our operator fails to satisfy the Fredholm property, the
solvability conditions can be formulated similarly. But this similarity is only formal since the range of the operator is not closed. The orthogonality
conditions~\eqref{ocb}  are with respect to the standard Fourier harmonics,
which solve the homogeneous adjoint problem for~\eqref{eq1}  when the
scalar potential function is trivial. They
belong to $L^{\infty}({\mathbb R}^{d})$ but they are not square-integrable.

We mention  that elliptic equations involving
non-Fredholm operators have been treated actively in recent years.
Approaches in weighted Sobolev and H\"older spaces have been developed in~\cite{Amrouche1997,Amrouche2008,Bolley1993,Bolley2001,B88}. In particular,
when $a=0$, the operator $A$ is Fredholm in certain properly chosen
weighted
spaces (see~\cite{Amrouche1997,Amrouche2008,Bolley1993,Bolley2001,B88}). However, the case where $a \neq 0$ is significantly
different and the approaches developed in these works are inapplicable.
The non-Fredholm Schr\"odinger type operators
were studied with the methods of the spectral and the
scattering theory in~\cite{EV21,V2011,VV08,VV19}.
Fredholm structures, topological invariants, and applications were considered
in~\cite{E09}. The works~\cite{GS05} and~\cite{RS01} are devoted to the
Fredholm and properness properties of   quasilinear elliptic
systems of the second order and of the operators of this kind on
${\mathbb R}^{N}$. The exponential decay and Fredholm properties in  
second-order quasilinear elliptic systems of equations were discussed in~\cite{GS10}.
The nonlinear non-Fredholm elliptic problems were treated in~\cite{EV20,VV14,VV21}. Significant applications to the theory of
reaction-diffusion
equations were developed in~\cite{DMV05,DMV08}.

The logarithmic Laplacian $\hbox{ln}(-\Delta)$ is the operator with Fourier
symbol $2\ln |p|$. It appears as the formal derivative
$\partial_{s}|_{s=0}(-\Delta)^{s}$ of fractional Laplacians at $s=0$.
The operator $(-\Delta)^{s}$ has been extensively  used, for instance, in the studies of  anomalous diffusion problems (to name a few, see~\cite{MK00,VV19,VV21} and  references therein). Spectral properties of the
logarithmic Laplacian in an open set of finite measure with Dirichlet boundary
conditions were  considered in~\cite{LW21} (see also ~\cite{CW19}).
The studies of $\ln (-\Delta)$ are important to understand the asymptotic spectral properties of the family of fractional Laplacians
in the limit $s\to 0^{+}$. In~\cite{JSW20}, it was demonstrated that such
an operator allows us to characterize the $s$-dependence of  solutions to
fractional Poisson problems for the full range of exponents $s\in (0, 1)$.

The operator contained in the left side of problem~\eqref{h}  is given by
\begin{equation}
\label{lab}
L_{a, b}:=\frac{1}{2}\ln\Big(-\frac{d^{2}}{dx^{2}}\Big)-b\frac{d}{dx}-a,
\quad a,b\in {\mathbb R}, \quad b\neq 0
\end{equation}
and is considered on $L^{2}({\mathbb R})$. By virtue of the standard Fourier
transform, it can be trivially obtained that the essential spectrum of
(\ref{lab}) equals to
\begin{equation}
\label{ess}
\lambda_{a,b}(p)=\ln \Big(\frac{|p|}{e^{a}}\Big)-ibp, \quad
a,b\in {\mathbb R}, \quad b\neq 0.
\end{equation}
Obviously, for $p\in {\mathbb R}$,
\begin{equation}
\label{lbab}
|\lambda_{a,b}(p)|=\sqrt{\ln^{2}\Big(\frac{|p|}{e^{a}}\Big)+b^{2}p^{2}}\geq
C_{a,b}>0,
\end{equation}
where $C_{a,b}$ is a constant. Note that as distinct from the situation without
the transport term considered in~\cite{EV23}, operator~\eqref{lab}  satisfies
the Fredholm property. Solvability of certain nonhomogeneous linear equations
containing the logarithmic Schr\"odinger operators in higher dimensions was
discussed in~\cite{EV231}. The article~\cite{ZKR23} deals with the symmetry
of positive solutions for Lane-Emden systems involving the logarithmic
Laplacian.

Let us set $K(x) = \varepsilon {\cal K}(x)$ with $\varepsilon \geq 0$.
In the applications to biology, such small and nonnegative parameter
$\varepsilon$ has the meaning
that the integral production term is small with respect to the others, that is, the frequency of large mutations is small enough.
Our first assumption is as follows.

\noindent {\bf Assumption 1.}  {\it The constants $a$, $b\in {\mathbb R}$ with $b\neq 0$.
The function $f(x): {\mathbb R}\to {\mathbb R}$ is  nontrivial such that
$f(x)\in L^{2}({\mathbb R})$. We also assume that
${\cal K}(x): {\mathbb R}\to {\mathbb R}$ does not vanish identically on
the real line and 
${\cal K}(x)\in L^{1}({\mathbb R})$.}

When the nonnegative parameter $\varepsilon=0$, we arrive at the
linear equation~\eqref{lp}.
By Lemma 5 below along with Assumption 1 and Remark 6,
problem~\eqref{lp}  admits a unique and nontrivial solution
\begin{equation}
\label{u0}
u_{0}(x)\in L^{2}({\mathbb R}).
\end{equation}
We look for the resulting solution of the nonlinear equation~\eqref{h}  as
\begin{equation}
\label{r}
u(x)=u_{0}(x)+u_{p}(x).
\end{equation}
Evidently, we easily derive the perturbed equation
\begin{equation}
\label{pert}
\Big[\frac{1}{2}\ln \Big(-\frac{d^{2}}{d x^{2}}\Big)\Big]u_{p}-
b\frac{du_{p}}{dx}-au_{p}
=\varepsilon \int_{-\infty}^{\infty}
{\cal K}(x-y)g(u_{0}(y)+u_{p}(y))dy.
\end{equation}
For the technical purposes, we will use a closed ball in our function space,
namely,
\begin{equation}
\label{b}
B_{\rho}:=\{u(x)\in L^{2}({\mathbb R}) \ | \ \|u\|_{L^{2}({\mathbb R})}\leq
\rho \}, \quad 0<\rho\leq 1.
\end{equation}
Let us seek  solutiona of equation~\eqref{pert} as  fixed points of the
auxiliary nonlinear problem,
\begin{equation}
\label{aux}
\Big[\frac{1}{2}\ln \Big(-\frac{d^{2}}{d x^{2}}\Big)\Big]u-
b\frac{du}{dx}-au
=\varepsilon \int_{-\infty}^{\infty}
{\cal K}(x-y)g(u_{0}(y)+v(y))dy,
\end{equation}
in the ball $B_{\rho}$ defined by~\eqref{b}. The fixed point technique was applied in~\cite{CV21} to study the
persistence of pulses for a class of reaction-diffusion equations. For a given function $v(y)$, \eqref{aux} is an equation to be solved for $u(x)$. 

The left side of~\eqref{aux}  contains the Fredholm operator introduced in~\eqref{lab}. The solvability of the linear and nonhomogeneous problem
involving such an operator will be discussed in the final section of this paper. We point out that, distinct from the present situation, the equations discussed in~\cite{EV20,VV14} involved
operators without the Fredholm property and orthogonality
conditions are required to solve them.

Now, we introduce the closed ball in the function space,
\begin{equation}
\label{M}
D_{M}:=\{g'(z)\in L^{\infty}({\mathbb R}) \ | \
\|g'(z)\|_{L^{\infty}({\mathbb R})}\leq M \}, \quad M>0
\end{equation}
and impose the following conditions on the nonlinear part of problem~\eqref{h}.

\noindent
{\bf Assumption 2.} {\it Let $g(z): {\mathbb R}\to {\mathbb R}$ such that
$g(0)=0$. Additionally, $g'(z)\in D_{M}$ and $g(z)$
does not vanish identically on the real line.}

We introduce the operator $t_g$ by $u = t_{g}v$, where $u$ is a
solution of equation~\eqref{aux}. As we will see, for $v\in B_{\rho}$, \eqref{aux} has only one solution $u\in B_{\rho}$.  Our first main statement is as
follows.

\noindent {\bf Theorem 3.} {\it Let Assumptions 1 and 2 hold and
\begin{equation}
\label{eps}
0<\varepsilon\leq \frac{\rho C_{a, b}}{M\|{\cal K}\|_{L^{1}({\mathbb R})}
(\|u_{0}\|_{L^{2}({\mathbb R})}+1)}.
\end{equation}
Then the map $t_{g}: B_{\rho}\to B_{\rho}$ associated with problem~\eqref{aux}
is a strict contraction.
The unique fixed point $u_{p}(x)$ of this map $t_{g}$ is the only solution of
equation~\eqref{pert}  in $B_{\rho}$.}

Clearly, the resulting solution of problem~\eqref{h}  given by~\eqref{r}  will
be nontrivial because the source term $f(x)$ does not vanish identically
on ${\mathbb R}$ and $g(0)=0$ as assumed.

Our second main result is about the continuity of the cumulative solution
of equation~\eqref{h}  given by~\eqref{r}  with respect to
the derivative of the nonlinear function $g$.
Let us define the following positive  auxiliary quantity
\begin{equation}
\label{sig}
\sigma:=\frac{\varepsilon M\|{\cal K}\|_{L^{1}({\mathbb R})}}{C_{a, b}}.
\end{equation}

\noindent
{\bf Theorem 4.} {\it Suppose that the assumptions of Theorem 3 hold.
Let $u_{p,j}(x)$ be the unique fixed point of the map
$t_{g_{j}}: B_{\rho}\to B_{\rho}$ for $g_{j}$, $j=1$, $2$ and the resulting solution of
problem~\eqref{h}  with $g(z)=g_{j}(z)$ be given by
\begin{equation}
\label{cum}
u_{j}(x):=u_{0}(x)+u_{p,j}(x).
\end{equation}
Then
\begin{equation}
\begin{aligned}
&\|u_{1}(x)-u_{2}(x)\|_{L^{2}({\mathbb R})}
\\
&\phantom {u_1(x)}\leq
\frac{\varepsilon}{1-
\sigma}\frac{\|{\cal K}\|_{L^{1}({\mathbb R})}}{C_{a, b}}
(\|u_{0}\|_{L^{2}({\mathbb R})}+1)\|g_{1}'(z)-g_{2}'(z)\|_{L^{\infty}({\mathbb R})}.
\end{aligned}
\label{cont}
\end{equation}
}

We proceed to the proof of our first main result.

%%%%%%%%%%%%%%%%%%%%%%%%%%%%%%%%%%%%%%%%%%%%%%%%%%
\bigskip

\setcounter{section}{2}
\setcounter{equation}{0}

\centerline{\bf 2. The existence of the perturbed solution}

{\it Proof of Theorem 3.} Let us choose an arbitrary $v(x)\in B_{\rho}$ and
denote the term contained in the integral expression in the right side of
equation~\eqref{aux}  as
\begin{equation}
\label{G}
G(x):=g(u_{0}(x)+v(x)).
\end{equation}
Throughout the work we will use the standard Fourier transform
\begin{equation}
\label{f}
\widehat{\phi}(p):=\frac{1}{\sqrt{2\pi}}\int_{-\infty}^{\infty}\phi(x)
e^{-ipx}dx.
\end{equation}
Obviously, the inequality
\begin{equation}
\label{fub}
\|\widehat{\phi}(p)\|_{L^{\infty}({\mathbb R})}\leq \frac{1}{\sqrt{2\pi}}
\|\phi(x)\|_{L^{1}({\mathbb R})}
\end{equation}
holds.
We apply (\ref{f}) to both sides of problem (\ref{aux}) and arrive at
\begin{equation}
\label{uhpKG}
\widehat{u}(p)=\varepsilon \sqrt{2\pi}
\frac{\widehat{\cal K}(p)\widehat{G}(p)}{\ln \Big(\frac{|p|}{e^{a}}\Big)
-ibp}.
\end{equation}
By means of inequalities (\ref{lbab}) and (\ref{fub}), we easily derive
\begin{equation}
\label{uhub}
|\widehat{u}(p)|\leq \varepsilon \frac{\|{\cal K}\|_{L^{1}({\mathbb R})}
|\widehat{G}(p)|}{C_{a, b}} 
\end{equation}
and then
\begin{equation}
\label{ul2ube}
\|u\|_{L^{2}({\mathbb R})}\leq \varepsilon
\frac{\|{\cal K}\|_{L^{1}({\mathbb R})}}{C_{a, b}}\|G\|_{L^{2}({\mathbb R})}.
\end{equation}
Clearly,
\begin{equation}
\label{Gx}
G(x)=\int_{0}^{u_{0}(x)+v(x)}g'(z)dz.
\end{equation}
Thus
\begin{equation}
\label{Gaub}
|G(x)|\leq  M|u_{0}(x)+v(x)|.
\end{equation}
It follows  that, for $v(x)\in B_{\rho}$, we have
\begin{equation}
\label{Gl2b}
\|G(x)\|_{L^{2}({\mathbb R})}\leq M
(\|u_{0}\|_{L^{2}({\mathbb R})}+1).
\end{equation}
Let us combine estimates (\ref{ul2ube}) and (\ref{Gl2b}) to obtain
\begin{equation}
\label{ul2rho}
\|u\|_{L^{2}({\mathbb R})}\leq \varepsilon \frac{\|{\cal K}\|_{L^{1}({\mathbb R})}}
{C_{a, b}}M(\|u_{0}\|_{L^{2}({\mathbb R})}+1)\leq \rho
\end{equation}
for all the values of the parameter $\varepsilon$  satisfying~\eqref{eps}, which means that   $u(x)\in B_{\rho}$ as well.

We suppose that for a certain $v(x)\in B_{\rho}$ there exist two solutions
$u_{1, 2}(x)\in B_{\rho}$ of problem~\eqref{aux}. Then the difference function
$w(x):=u_{1}(x)-u_{2}(x)\in L^{2}({\mathbb R})$ is a solution of the
homogeneous equation~\eqref{lnw}. Thus $w(x)$ is trivial on ${\mathbb R}$
as discussed with similar  arguments for Lemma 5 below.

Therefore, problem~\eqref{aux}  defines a map
$t_{g}: B_{\rho}\to B_{\rho}$ for all the values of $\varepsilon$ satisfying~\eqref{eps}.
Our goal is to demonstrate that this map is a strict contraction in the following.

Let us choose
arbitrary $v_{1,2}(x)\in B_{\rho}$. By virtue of the argument above
$u_{1,2}:=t_{g}v_{1,2}\in B_{\rho}$ as well when $\varepsilon$ satisfies~\eqref{eps}.
From~\eqref{aux}  we easily deduce that
\begin{equation}
\label{aux1}
\left[\frac{1}{2}\ln \left(-\frac{d^{2}}{dx^{2}}\right)\right] u_{1}-
b\frac{du_{1}}{dx}-au_{1} = 
\varepsilon \int_{-\infty}^{\infty}{\cal K}(x-y)g(u_{0}(y)+v_{1}(y))dy,
\end{equation}
\begin{equation}
\label{aux2}
\left[\frac{1}{2}\ln \left(-\frac{d^{2}}{dx^{2}}\right)\right] u_{2}-
b\frac{du_{2}}{dx}-au_{2} = 
\varepsilon \int_{-\infty}^{\infty}{\cal K}(x-y)g(u_{0}(y)+v_{2}(y))dy.
\end{equation}
Let us introduce
\begin{equation}
\label{G1G2}
G_{1}(x):=g(u_{0}(x)+v_{1}(x)), \quad G_{2}(x):=g(u_{0}(x)+v_{2}(x))
\end{equation}
and apply the standard Fourier transform (\ref{f}) to both sides of
equations (\ref{aux1}) and (\ref{aux2}). This gives us
\begin{equation}
\label{u12hpKG}
\widehat{u_{1}}(p)=\varepsilon \sqrt{2\pi}
\frac{\widehat{\cal K}(p)\widehat{G_{1}}(p)}{\ln \Big(\frac{|p|}{e^{a}}
\Big)-ibp}, \quad
\widehat{u_{2}}(p)=\varepsilon \sqrt{2\pi}
\frac{\widehat{\cal K}(p)\widehat{G_{2}}(p)}{\ln \Big(\frac{|p|}{e^{a}}
\Big)-ibp}.
\end{equation}
By means of bounds (\ref{lbab}) and (\ref{fub}), we arrive at
\begin{equation}
\label{u12hub}
|\widehat{u_{1}}(p)-\widehat{u_{2}}(p)|\leq \varepsilon
\|{\cal K}\|_{L^{1}({\mathbb R})}
\frac{|\widehat{G_{1}}(p)-\widehat{G_{2}}(p)|}{C_{a, b}}.
\end{equation}
Thus
\begin{equation}
\label{u12ubn}
\|u_{1}-u_{2}\|_{L^{2}({\mathbb R})}\leq \varepsilon
\frac{\|{\cal K}\|_{L^{1}({\mathbb R})}}{C_{a, b}}\|G_{1}-G_{2}\|_{L^{2}({\mathbb R})}.
\end{equation}
Note that
\begin{equation}
\label{G12i}
G_{1}(x)-G_{2}(x)=\int_{u_{0}(x)+v_{2}(x)}^{u_{0}(x)+v_{1}(x)}g'(z)dz.
\end{equation}
This enables us to estimate the norm as
\begin{equation}
\label{G12nub}
\|G_{1}-G_{2}\|_{L^{2}({\mathbb R})}\leq M\|v_{1}-v_{2}\|_{L^{2}({\mathbb R})}.
\end{equation}
We use inequalities (\ref{u12ubn}) and (\ref{G12nub}) along with
definition (\ref{sig}) to derive
\begin{equation}
\label{contr}
\|t_{g}v_{1}-t_{g}v_{2}\|_{L^{2}({\mathbb R})}\leq \sigma \|v_{1}-v_{2}\|_
{L^{2}({\mathbb R})}.
\end{equation}
By virtue of condition (\ref{eps}) along with Remark 6 below, we have
\begin{equation}
\label{sig1}
\sigma<1.
\end{equation}
This implies that the map $t_{g}: B_{\rho}\to B_{\rho}$ defined by equation~\eqref{aux}  is a strict contraction for all the values of $\varepsilon$
satisfying~\eqref{eps}. Its unique fixed point $u_{p}(x)$ is the only
solution of problem~\eqref{pert} in the ball $B_{\rho}$. 
The resulting
$u(x)\in L^{2}(\mathbb R)$ given by formula~\eqref{r}  solves
equation~\eqref{h}. Note that by virtue of~\eqref{ul2rho}, we have that $\|u_{p}(x)\|_{L^{2}(\mathbb R)}\to 0$ as
$\varepsilon\to 0$. \hfill\lanbox

Now we turn our attention to  establish  the second main
result of this paper.
 
\bigskip

%%%%%%%%%%%%%%%%%%%%%%%%%%%%%%%%%%%%%%%%%%%%%%%%%%

\setcounter{section}{3}
\setcounter{equation}{0}

\centerline{\bf 3. The continuity of the resulting solution}

\bigskip

\noindent
{\it Proof of Theorem 4.} Evidently, for all the values of
$\varepsilon$ which satisfy~\eqref{eps}, we have
\begin{equation}
\label{up12}
u_{p,1}=t_{g_{1}}u_{p,1}, \quad u_{p,2}=t_{g_{2}}u_{p,2}.
\end{equation}
Noting 
\begin{equation}
\label{up12m}
u_{p,1}-u_{p,2}=t_{g_{1}}u_{p,1}-t_{g_{1}}u_{p,2}+t_{g_{1}}u_{p,2}-
t_{g_{2}}u_{p,2},
\end{equation}
we get
\begin{equation}
\label{up12ni}
\|u_{p,1}-u_{p,2}\|_{L^{2}({\mathbb R})}\leq\|t_{g_{1}}u_{p,1}-t_{g_{1}}
u_{p,2}\|_{L^{2}({\mathbb R})}+\|t_{g_{1}}u_{p,2}-t_{g_{2}}u_{p,2}\|_
{L^{2}({\mathbb R})}.
\end{equation}
By means of estimate (\ref{contr}),
\begin{equation}
\label{tg1up12}
\|t_{g_{1}}u_{p,1}-t_{g_{1}}u_{p,2}\|_{L^{2}({\mathbb R})}\leq
\sigma\|u_{p,1}-u_{p,2}\|_{L^{2}({\mathbb R})}
\end{equation}
and inequality (\ref{sig1}) holds.
Then we arrive at
\begin{equation}
\label{sigma}
(1-\sigma)\|u_{p,1}-u_{p,2}\|_{L^{2}({\mathbb R})}\leq
\|t_{g_{1}}u_{p,2}-t_{g_{2}}u_{p,2}\|_{L^{2}({\mathbb R})}.
\end{equation}
Let $\xi(x):=t_{g_{1}}u_{p,2}$. Clearly,
\begin{equation}
\label{12}
\Big[\frac{1}{2}\ln \Big(-\frac{d^{2}}{dx^{2}}\Big)\Big]\xi-
b\frac{d\xi}{dx}-a\xi=\varepsilon
\int_{-\infty}^{\infty}{\cal K}(x-y)g_{1}(u_{0}(y)+u_{p,2}(y))dy,
\end{equation}
\begin{equation}
\label{22}
\Big[\frac{1}{2}\ln \Big(-\frac{d^{2}}{dx^{2}}\Big)\Big]u_{p,2}-
b\frac{du_{p,2}}{dx}-au_{p,2}=\varepsilon
\int_{-\infty}^{\infty}{\cal K}(x-y)g_{2}(u_{0}(y)+u_{p,2}(y))dy.
\end{equation}
Denote
\begin{equation}
\label{G12G22}
G_{1,2}(x):=g_{1}(u_{0}(x)+u_{p,2}(x)), \quad G_{2,2}(x):=g_{2}(u_{0}(x)+u_{p,2}(x)).
\end{equation}
By applying the standard Fourier transform (\ref{f}) to both sides of
equations (\ref{12}) and (\ref{22}), we obtain
\begin{equation}
\label{xiup2h}
\widehat{\xi}(p)=\varepsilon \sqrt{2 \pi}\frac{\widehat{\cal K}(p)
\widehat{G_{1,2}}(p)}{\ln \Big(\frac{|p|}{e^{a}}\Big)-ibp}, \quad
\widehat{u_{p,2}}(p)=\varepsilon \sqrt{2 \pi}\frac{\widehat{\cal K}(p)
\widehat{G_{2,2}}(p)}{\ln \Big(\frac{|p|}{e^{a}}\Big)-ibp}.
\end{equation}
With the assistance of~\eqref{lbab}  and~\eqref{fub}, we get
\begin{equation}
\label{xiup2hub}
|\widehat{\xi}(p)-\widehat{u_{p,2}}(p)|\leq \varepsilon
\frac{\|{\cal K}\|_{L^{1}({\mathbb R})}}{C_{a, b}}
|\widehat{G_{1,2}}(p)-\widehat{G_{2,2}}(p)|.
\end{equation}
Then
\begin{equation}
\label{xiup2hubn}
\|\xi-u_{p,2}\|_{L^{2}({\mathbb R})}\leq \varepsilon
\frac{\|{\cal K}\|_{L^{1}({\mathbb R})}}{C_{a, b}}\|G_{1,2}-G_{2,2}\|_{L^{2}({\mathbb R})}.
\end{equation}
Obviously,
$$
G_{1,2}(x)-G_{2,2}(x)=\int_{0}^{u_{0}(x)+u_{p,2}(x)}[g_{1}'(z)-g_{2}'(z)]dz.
$$
This allows us to derive an upper bound on the norm as
\begin{equation}
\label{g1222l2n}
\|G_{1,2}-G_{2,2}\|_{L^{2}({\mathbb R})}\leq
\|g_{1}'(z)-g_{2}'(z)\|_{L^{\infty}({\mathbb R})}(\|u_{0}\|_{L^{2}({\mathbb R})}+1).
\end{equation}
By means of inequalities (\ref{sigma}), (\ref{xiup2hubn}) and (\ref{g1222l2n}),
the norm
$\|u_{p,1}-u_{p,2}\|_{L^{2}({\mathbb R})}$ can be estimated from above by
\begin{equation}
\label{up1p2h1}
\frac{\varepsilon}{1-\sigma}\frac{\|{\cal K}\|_{L^{1}({\mathbb R})}}{C_{a, b}}
(\|u_{0}\|_{L^{2}({\mathbb R})}+1)\|g_{1}'(z)-g_{2}'(z)\|_{L^{\infty}({\mathbb R})}.
\end{equation}
Thus~\eqref{cont} follows easily from~\eqref{cum} and~\eqref{up1p2h1}.       \hfill\lanbox

\bigskip

%%%%%%%%%%%%%%%%%%%%%%%%%%%%%%%%%%%%%%%%%%%%%%%%%%

\setcounter{section}{4}
\setcounter{equation}{0}

\centerline{\bf 4. Auxiliary results}

We address the solvability of the linear problem containing the
logarithmic Laplacian, a drift term, and a square-integrable right
side,
\begin{equation}
\label{lp}
\Big[\frac{1}{2}\ln \Big(-\frac{d^{2}}{dx^{2}}\Big)\Big]u-
b\frac{du}{dx}-au=f(x), \quad
x\in {\mathbb R},
\end{equation}
where $a, b\in {\mathbb R}, \ b\neq 0$ are constants.

\noindent
{\bf Lemma 5.} {\it  Let  $f(x):{\mathbb R}\to {\mathbb R}$ be nontrivial and
$f(x)\in L^{2} ({\mathbb R})$. Then equation (\ref{lp}) possesses a unique
solution $u_{0}(x)\in L^{2} ({\mathbb R})$.}

\bigskip

\noindent
{\it Proof.} To demonstrate the uniqueness of solutions for
problem (\ref{lp}), we suppose that it has two solutions
$u_{1}(x)$, $u_{2}(x)\in L^{2}({\mathbb R})$.
Clearly, the difference $w(x):=u_{1}(x)-u_{2}(x)\in L^{2}({\mathbb R})$ as
well, which   solves the homogeneous equation,
\begin{equation}
\label{lnw}
\Big[\frac{1}{2}\ln \Big(-\frac{d^{2}}{dx^{2}}\Big)\Big]w-b\frac{dw}{dx}
-aw=0.
\end{equation}
Since the operator $L_{a, b}$ on $L^{2}({\mathbb R})$
defined by~\eqref{lab}  has only the essential spectrum and no
nontrivial zero modes (see~\eqref{ess} and~\eqref{lbab}),
the function $w(x)$ vanishes a.e. on the real line.

We apply the standard Fourier transform~\eqref{f}  to both sides of
problem~\eqref{lp} to get
\begin{equation}
\label{upfp0}
\widehat{u}(p)=\frac{\widehat{f}(p)}
{\ln \Big(\frac{|p|}{e^{a}}\Big)-ibp}, \quad p\in {\mathbb R}.
\end{equation}
By virtue of (\ref{lbab}), we can derive that
\[
|\widehat{u}(p)|\leq \frac{|\widehat{f}(p)|}{C_{a, b}}\in
L^{2}({\mathbb R})
\]
via our assumption. Thus $u(x)\in L^{2}({\mathbb R})$. \hfill\lanbox

\bigskip

\noindent
{\bf Remark 6.} {\it Under the conditions of Lemma 5, the unique square-integrable solution $u_{0}(x)$ of equation (\ref{lp}) does not
vanish identically on ${\mathbb R}$ since the source term $f(x)$ is nontrivial.}

\bigskip

%%%%%%%%%%%%%%%%%%%%%%%%%%%%%%%%%%%%%%%%%%%%%%%%%%%%%%%%%%%%%%%%

\section*{Acknowledgements}

Y.C. was supported by  NSERC of Canada (Nos.\ RGPIN-2019-05892, RGPIN-2024-05593).

\end{document}